\documentclass [12pt,a4paper]{amsart}
\usepackage{}
\usepackage{stmaryrd}
\usepackage{textcomp}
\usepackage{amsmath}
\usepackage{amssymb}
\usepackage{cases}
\usepackage{amscd}
\usepackage{epsfig}
\usepackage{graphicx}
\usepackage{verbatim}
\usepackage{amssymb}
\usepackage{amsfonts}
\usepackage{amsbsy}
\usepackage{amsthm}
\usepackage{amstext}
\usepackage{amsopn}
\usepackage[dvips]{color}
\usepackage{indentfirst}
\usepackage{mathrsfs}
\usepackage{enumerate}
\newtheorem{thm}{Theorem}[section]

\newtheorem{lem}[thm]{Lemma}
\newtheorem{prop}[thm]{Proposition}
\newtheorem{rem}[thm]{Remark}

\newcommand{\PU}{\textrm{PU}}
\newcommand{\SU}{\textrm{SU}}
\newcommand{\HC}{{\textbf{H}_{\mathbb{C}}^2}}

\renewcommand{\Re}{\,{\operatorname{Re}\,}}
\newcommand{\tr}{\textrm{tr}}

\newcommand{\Det}{\textrm{Det}}

\begin{document}

\author{John R. Parker and Li-jie Sun}
\address{
University of Durham\\
South Road\\
Durham DH1 3LE UK}
\email{j.r.parker@durham.ac.uk}
\address{
Tokyo Institute of Technology\\
2-12-1 O-okayama Meguro-ku\\
Tokyo 152-8552 Japan}
\email{lijie.alice.sun@gmail.com}
\title{complex hyperbolic triangle groups with 2-fold symmetry}

\begin{abstract}
In this paper we will consider the 2-fold symmetric complex hyperbolic triangle groups generated by 
three complex reflections through angle $2\pi/p$ with $p\geqslant2$. We will mainly concentrate on the groups where 
some elements are elliptic of finite order. Then we will classify all such groups which are candidates for being 
discrete. There are only 4 types.
\end{abstract}
\maketitle
\section{ Introduction}

A complex hyperbolic triangle is a triple $(C_{1}, C_{2}, C_{3})$ of complex geodesics in $\HC$. If each pair 
of complex geodesics intersects in $\HC\cup\partial\HC$ and the angles between $C_{k-1}$ and $C_{k}$ for 
$k=1, 2, 3$ (the indices are taken mod 3) are $\pi/p_{1}, \pi/p_{2}, \pi/p_{3}$, where 
$p_{1}, p_{2}, p_{3} \in\mathbb{N}\cup\{\infty\}$, we call the triangle $(C_{1}, C_{2}, C_{3})$ a 
$(p_{1}, p_{2}, p_{3})$-$triangle$. The intersection points of pairs of complex geodesics are called the vertices 
of the complex hyperbolic triangle. A group $\Gamma$ is called a $(p_{1}, p_{2}, p_{3})-$\emph{triangle group}, 
if $\Gamma$ is generated by three complex reflections $R_{1}, R_{2}, R_{3}$ fixing sides $C_{1},C_{2},C_{3}$ of 
$(p_{1}, p_{2}, p_{3})$-$triangle$. Note that a complex reflection may have order greater than $2$. In what
follows we suppose that $R_1$, $R_2$ and $R_3$ all have order $p\in{\mathbb Z}$ with $p\geqslant 2$.

\par Any two real hyperbolic triangle groups with the same intersection angles are conjugate in 
$\textrm{Isom}^{+}(\textbf{H}^3)$, which is the orientation preserving isometry group of 3-dimensional 
hyperbolic space. If we consider the groups in $\PU(2, 1)=\textrm{Aut}(\HC)$, we will get the nontrivial 
deformations. The deformation theory of complex hyperbolic triangle groups was begun in \cite{GP1} in which 
they investigated $\Gamma$ of type $(\infty,\infty,\infty)$ with $p=2$ (complex hyperbolic ideal triangle groups).  
Since then, there have been many developments referring to other types, such as  \cite{Wy, KPT2, Sun} 
among which they mainly gave the necessary conditions of $\Gamma$ to be discrete. Especially Parker and 
Paupert in \cite{Par3} and \cite{PP} investigated the equilateral triangle group generated by three complex 
reflections with finite order. These include \emph{Deraux's lattice, Livn\'e's lattices, Mostow's lattices}. 
Our starting point is a result given by Thompson \cite{Thom} where he investigated the non-equilateral triangle 
groups generated by three complex involutions (that is the order of the reflections is $p=2$). He obtained 
his result using a computer search. Using \cite{PP} we see that Thompson's results apply to groups 
with $p>2$ as well. In what follows we will give the specific case about the triangles group with 
2-fold symmetry and we give a rigorous proof. 

We will restrict to the complex hyperbolic triangle groups generated by three complex reflections with finite
order $p\ge 2$. Suppose that the polar vector of a complex geodesic $C_{1}$ is $\bf v_{1}$ (see Section~2 
for a more precise explanation). We consider the complex reflection $R_{1}$ in the complex geodesic 
$C_{1}$. This map sends $\bf v_{1}$ to $e^{i\phi}\bf v_{1}$ and acts 
as the identity on the orthogonal complement of ${\bf v}_1$, that is on vectors that project to $C_1$. 
We will always restrict to the case where $\phi=2\pi/p$ and so $R_1$ has order $p\geqslant2$. 
\begin{equation}\label{ref:l}
R_{1}({\bf z})={\bf z}+(e^{i\phi}-1)
\frac{\langle {\bf z}, {\bf v_1}\rangle}{\langle {\bf v_1}, {\bf v_1}\rangle}{\bf v_1}.
\end{equation}
In order to convert $R_1$ into a matrix with determinant 1, we need to multiply 
the expression in \eqref{ref:l} by $e^{-i\phi/3}$.
The ambiguity involved in this choice is precisely the ambiguity involved in lifting an isometry in
$\PU(2,1)$ to a matrix in $\SU(2,1)$. 

Here we recall the terminology for \emph{braid relations} between group elements (see Section 2.2 of Mostow 
\cite{Mos}). Let $G$ be a group and $a,\, b\in G$. Then $a$ and $b$ satisfy a braid relation of length 
$l\in \mathbb{Z}_{+}$ if 
\[
(ab)^{l/2}=(ba)^{l/2},
\]
where powers means that the corresponding alternating product of $a$ and $b$ should have $l$ factors. For 
example, $(ab)^{3/2}=aba$, $(ba)^2=baba$. We denote the \emph{braid length} $l$ by $br (a, b)$ to be the 
minimum length of a braid relation satisfied by $a$ and $b$.

We define the $(l_{1}, l_{2}, l_{3}; l_{4})$-triangle groups to be the triangle groups 
with the following braid relations: 
$$
\begin{array}{ll}
br(R_2, R_3)=l_1,\quad & br(R_1, R_3)=l_2, \\
br(R_1, R_2)=l_3,\quad & br(R_1, R_3^{-1}R_2R_3)=l_4,
\end{array}
$$
where $R_{j}$ is of order $p$.

In this paper we aim to list the candidates of discrete triangle groups generated by $R_{1}, R_{2}, R_{3}$ with 
$l_{1}=l_{2}$ and $l_{3}=l_{4}$ as stated in Theorem \ref{thm:main}. 

\section{The parameter space, traces and main result}

\par Firstly we recall some fundamentals about complex hyperbolic 2-space. Please refer to \cite{Gol, Par4} for 
more details about the complex hyperbolic space. Let $\mathbb{C}^{2, 1}$ denote the vector space 
$\mathbb{C}^{3}$ equipped with the Hermitian form
$$\langle {\bf z}, {\bf w} \rangle=z_{1}\overline{w}_{1}+z_{2}\overline{w}_{2}-z_{3}\overline{w}_{3}$$
of signature (2,1), where ${\bf z}=[z_{1}, z_{2}, z_{3}]^{t}$ and ${\bf w}=[w_{1}, w_{2}, w_{3}]^{t}.$ The 
Hermitian form divides 
$\mathbb{C}^{2,1}$ into three parts $V_{-}, V_{0}$ and $V_{+}$, which are
\begin{align*}
&V_{-}=\{{\bf z}\in\mathbb{C}^{2,1}| \langle {\bf z}, {\bf z}\rangle<0\},\\
&V_{0}=\{{\bf z}\in\mathbb{C}^{2,1}|\langle{\bf z}, {\bf z}\rangle=0\},\\ 
&V_{+}=\{{\bf z}\in\mathbb{C}^{2,1}| \langle {\bf z}, {\bf z}\rangle>0\}.
\end{align*}

We denote by $\mathbb{CP}^{2}$ the complex projectivisation of $\mathbb{C}^{2,1}$ and by $\mathbb{P}:$ 
$\mathbb{C}^{2,1}\setminus \{0\}\rightarrow\mathbb{CP}^{2}$ the natural projectivisation map. The 
\emph{complex hyperbolic 2-space} $\HC$ is defined as $\mathbb{P}(V_{-})$. It is called the standard projective 
model of the complex hyperbolic space. Correspondingly the boundary of $\HC$ is 
$\partial \HC=\mathbb{P}({V_{0}\setminus \{0\}})$. One can also consider the \emph{unit ball model} whose 
boundary is the sphere $\mathbb{S}^3$ by taking $z_{3}=1,$ which can be simply written as 
$\{(z_1, z_2)\in\mathbb{C}^2: |z_1|^2+|z_2|^2<1\}.$ 

\par The complex hyperbolic plane $\HC$ is a $K\ddot{a}hler$ manifold of constant holomorphic sectional curvature $-1$. The holomorphic automorphism group of $\HC$ is the projectivisation $\PU(2,1)$ of the group $\textrm{U}(2,1)$ of complex linear transformations on $\mathbb{C}^{2,1}$, which preserve the Hermitian form. Especially $\SU(2,1)$ is the subgroup of $\textrm{U}(2,1)$ with the determinant of each element being 1.

Let $x,  y\in \HC$ be points corresponding to vectors ${\bf x},{\bf y}\in \mathbb{C}^{2,1}\setminus \{0\}$. Then the \textit{Bergman metric} $\rho$ on $\HC$ is given by
$$\cosh^2\Big(\frac{\rho({ x}, {y})}{2}\Big)=\frac{\langle{\bf x},{\bf y}\rangle\langle {\bf y},{\bf x}\rangle}{\langle {\bf x},{\bf x}\rangle\langle {\bf y}, {\bf y}\rangle},$$
where ${\bf x}$, ${\bf y}{\in V_{-}}$ are the lifts of $ x,$ $y$ respectively. It is easy to check that this definition is independent of the choice of lifts.

    Given two points $ x$ and $ y$ in $\HC\cup \partial \HC$, with lifts ${\bf x}$ and ${\bf y}$ to $\mathbb{C}^{2,1}$ respectively, the complex span of ${\bf x}$ and ${\bf y}$ projects to a \emph{complex line} in $\mathbb{CP}^{2}$ passing through 
$x$ and $y$. The intersection of a complex line with $\HC$ will be called a \emph{complex geodesic} $C$, which can be uniquely determined by a positive vector ${\bf v}\in V_+,$ i.e. $C=\mathbb{P}(\{{\bf z}\in\mathbb{C}^{2,1}\setminus \{0\}|\langle {\bf z}, {\bf v}\rangle=0\})$. We call ${\bf v}$ a \emph{polar vector} to $C.$ As stated in Section 1, we will consider $(l_{1}, l_{2}, l_{3}; l_{4})$-triangle groups $\Gamma$ generated by three complex reflections (see (\ref{ref:l})) through angle $\phi$ in three complex geodesics.

Throughout this paper, we assume that $R_1,\, R_2,\, R_3$ are three complex reflections in complex geodesics $C_1,\, C_2,\, C_3$ respectively.
We parameterize the triangle groups generated by $R_{1},\, R_{2},\, R_{3}$ by three complex numbers $\rho$, $\sigma$ and $\tau$.
Up to the action of ${\rm PU}(2,1)$, we can parameterize the collection of three 
pairwise distinct complex lines in ${\bf H}^2_{\mathbb C}$ by four real parameters. 
The parameters we choose are $|\rho|$, $|\sigma|$, $|\tau|$ and 
${\rm arg}(\rho\sigma\tau)$. In particular, we can freely choose the argument
of two out of the three parameters.

Write $u=e^{i\phi/3}=e^{2\pi i/3p}$. The group $\Gamma$ has generators given by 
\begin{equation}\label{R123}
  R_1=\left(\begin{matrix}
  u^2 & \rho & -u\bar \tau\\
  0 & \bar u & 0\\
  0 & 0 & \bar u
  \end{matrix}\right),\\
  R_2=\left(\begin{matrix}
  \bar u & 0 & 0\\
  -u\bar \rho & u^2 & \sigma \\
  0 & 0 & \bar u
\end{matrix}\right),\\
  R_3=\left(\begin{matrix}
  \bar u & 0 & 0\\
  0 & \bar u & 0\\
  \tau & -u\bar \sigma & u^2
\end{matrix}\right)
\end{equation}
which preserve the Hermitian form
\begin{equation}\label{cal:H}
H=\left(\begin{matrix}
  \alpha & \beta_1 & \overline{\beta}_3\\
  \overline{\beta}_1&\alpha&\beta_2\\
  \beta_3&\overline{\beta}_2&\alpha
\end{matrix}\right)
\end{equation}
where $\alpha=\sqrt{2-u^3-\bar u^3}$, $\beta_1=-i\bar{u}^{1/2}\rho$,
$\beta_2=-i\bar{u}^{1/2}\sigma$, $\beta_3=-i\bar{u}^{1/2}\tau$ (note that here we take $\bar{u}^{1/2}=e^{-\pi i/3p}$).

This Hermitian form has signature $(2,1)$ if and only if ${\rm det}(H)<0$. That is,
\begin{eqnarray*}
0 & < & \alpha|\beta_1|^2+\alpha|\beta_2|^2+\alpha|\beta_3|^2
-\alpha^3-\beta_1\beta_2\beta_3-\bar{\beta}_1\bar{\beta}_2\bar{\beta}_3 \\
& = & \alpha^2|\rho|^2+\alpha^2|\sigma|^2+\alpha^2|\tau|^2
-\alpha^3-i\bar{u}^{3/2}\rho\sigma\tau+
i u^{3/2}\bar{\rho}\bar{\sigma}\bar{\tau}.
\end{eqnarray*}

In terms of these parameters
\begin{equation}\label{equa:tr}
\begin{array}{ll}
&{\rm tr}(R_1R_2)=u(2-|\rho|^2)+\bar{u}^2, \quad  
{\rm tr}(R_2R_3)=u(2-|\sigma|^2)+\bar{u}^2, \\
&{\rm tr}(R_1R_3)=u(2-|\tau|^2)+\bar{u}^2, \quad\\
&{\rm tr}(R_1R_3^{-1}R_2R_3)=u(2-|\sigma\tau-\bar{\rho}|^2)+\bar{u}^2.
\end{array}
\end{equation}

\begin{lem}\text{\cite[Corollary 2.5]{PP}}\label{lem:zeta}
If $|\rho|=2\cos\zeta$, then the three eigenvalues of $R_1R_2$ will be $\bar{u}^2, \, -ue^{2i\zeta} ,\, -ue^{-2i\zeta}$.
\end{lem}
\begin{proof} Each point on $C_1$ is a $\bar{u}=e^{-i\phi/3}$ eigenvector of $R_{1}$ and each point on 
$C_2$ is a $\bar{u}=e^{-i\phi/3}$ eigenvector of $R_{2}$, see (\ref{ref:l}). Therefore if 
${\bf z}\in C_1\cap C_2$, then we will get that
\[
R_1R_2({\bf z})=e^{-i\phi/3}R_{1}({\bf z})=e^{-2i\phi/3}\bf z
\]
Hence $\bf z$ is a $\bar{u}^2=e^{-2i\phi3}$ eigenvector of $R_1R_2$. Hence the sum of the other two
eigenvalues of $R_1R_2$ is $u(2-|\rho|^2)$.  By the assumption $|\rho|=2\cos\zeta$, we know that $R_1R_2$ is not loxodromic, see Section 6.2 in \cite{Gol}. Therefore each eigenvalue of  $R_1R_2$ is of modulus one. Then we can get that the 
three eigenvalues of $R_1R_2$ will be $\bar{u}^2, \, -ue^{2i\zeta} ,\, -ue^{-2i\zeta}$ from the form of 
$\tr(R_1R_2)$ in (\ref{equa:tr}). 
\end{proof}

\begin{rem}\label{rem:1}
If we suppose that $m\in\mathbb{N}$, $m \geqslant 2$, then $|\rho|=2\cos(\pi/m)$ if and only if $br(R_1, R_2)=m$, see Section 2.2 in \cite{Mos} for details or more precisely \text{\cite[prop 2.3]{DPP2}}. In particular if $R_1,\, R_2$ are complex involutions ($p=2$), then the order of $R_1R_2$ will be of $m$.
\end{rem}
Assume that 
\[
br (R_1,\, R_2) = br (R_1,\, R_3^{-1}R_2R_3),\qquad br (R_2,\, R_3) = br (R_1,\, R_3).
\]
From Remark \ref{rem:1} and (\ref{equa:tr}), our hypothesis on braiding implies that
$$
|\rho|=|\sigma\tau-\bar{\rho}|,\qquad |\sigma|=|\tau|.
$$
Since we are free to choose the argument of two of the three parameters, we impose
the condition that $\sigma$ and $\tau$ should be real and non-negative, which means that
${\rm Im}(\rho)={\rm Im}(\sigma\tau-\bar{\rho})$. So the condition
$|\rho|=|\sigma\tau-\bar{\rho}|$ becomes either $\sigma\tau=\rho+\bar{\rho}$ or $\sigma\tau=0$.
In the latter case the group is reducible, so we do not consider it. Hence we suppose
${\rm Re}(\rho)> 0$ and $\sigma=\tau=\sqrt{\rho+\bar{\rho}}$. 

We suppose that $|\rho|=2\cos(\pi/m)$ and $\sigma=\tau=2\cos(\pi/n)$, where $m,\, n\in\mathbb{N}$  and $m,\, n\geqslant 3$. Therefore the matrices in (\ref{R123}) become:
\begin{eqnarray}\label{mar:1}
  R_1  & = &   \left(\begin{matrix}
  u^2 & \rho & -u \sqrt{\rho+\bar{\rho}}  \\
  0 & \bar u & 0\\
  0 & 0 & \bar u
  \end{matrix}\right), \\
  \label{mar:2}
  R_2  & = &   \left(\begin{matrix}
  \bar u & 0 & 0\\
  -u\bar \rho & u^2 & \sqrt{\rho+\bar{\rho}} \\
  0 & 0 & \bar u
\end{matrix}\right),\\
\label{mar:3}
  R_3 & = &  \left(\begin{matrix}
  \bar u & 0 & 0\\
  0 & \bar u & 0\\
  \sqrt{\rho+\bar{\rho}} & -u\sqrt{\rho+\bar{\rho}} & u^2
\end{matrix}\right).
\end{eqnarray}
Furthermore, the Hermitian form $H$ (\ref{cal:H}) has signature $(2,1)$ if and only if
\begin{equation}\label{H:verify}
\begin{aligned}
0 & <  \alpha|\rho|^2+2\alpha(\rho+\bar{\rho})-\alpha^3-i\bar{u}^{3/2}(\rho^2+|\rho|^2)
+iu^{3/2}(\bar{\rho}^2+|\rho|^2) \\
& = 2\alpha(\rho+\bar{\rho})-\alpha^3-i\bar{u}^{3/2}\rho^2+i u^{3/2}\bar{\rho}^2.
\end{aligned}
\end{equation}
\begin{prop}\label{prop:rela}
Let
$$
S=\left(\begin{matrix}
\rho & u(1-\rho-\bar{\rho}) & u^2\sqrt{\rho+\bar{\rho}} \\
\bar{u} & 0 & 0 \\ 0 & \bar{u}\sqrt{\rho+\bar{\rho}} & -1 
\end{matrix}\right).
$$
Then
\begin{enumerate}
\item[(a)] $S^2=R_1R_2R_3$,
\item[(b)]
\begin{eqnarray*}
 SR_1S^{-1} &=& R_1R_2R_1^{-1}, \\ 
 SR_2S^{-1} &=& R_1R_3R_1R_3^{-1}R_1^{-1}, \\
 SR_3S^{-1} &=& R_1R_3R_1^{-1}.
 \end{eqnarray*}
In particular,
$$
S(R_2R_3)S^{-1}=R_1R_3,\quad S(R_1R_3^{-1}R_2R_3)S^{-1}=R_1R_2.
$$
\end{enumerate}

\end{prop}

\begin{proof}
The first part follows by a simple matrix calculation. The second part could be verified
similarly, but it is easy to observe that 
$$
S{\bf v}_1=R_1{\bf v}_2, \quad S{\bf v}_2= R_1R_3{\bf v}_1,\quad 
S{\bf v}_3=-uR_1{\bf v}_3,
$$
where ${\bf v}_{1}=[1, 0, 0]^{t},$ ${\bf v}_{2}=[0, 1, 0]^{t}$ and ${\bf v}_{3}=[0, 0, 1]^{t}$ are the polar 
vectors of $R_{1}, \, R_{2}, \, R_{3}$ respectively.
The result follows.
\end{proof}

In the following we will classify all discrete triangle groups generated by $R_{1}, R_{2}, R_{3}$ with the 2-fold symmetry given by 
$S$ satisfying the conditions (a) and (b) in Proposition \ref{prop:rela}.

\begin{thm}\label{thm:main}
Let $R_{1}, R_{2}, R_{3}$ be three complex reflections of order $p$ (with $p\geqslant2$) in $\SU(2,1)$ 
so that $R_i$ keeps a complex geodesic $C_i$
$(i=1,\, 2,\, 3)$ invariant. Assume that there is $S\in\SU(2, 1)$ such that
\[
SR_{1}S^{-1}=R_{1}R_{2}R_{1}^{-1},~ 
SR_{2}S^{-1}=R_{1}R_{3}R_{1}R_{3}^{-1}R_{1}^{-1}, 
SR_{3}S^{-1}=R_{1}R_{3}R_{1}^{-1},
\]
\[
S^{2}=R_{1}R_{2}R_{3}.
\]
Suppose that $br(R_1, R_3)=n$, $br(R_{1}, R_{2})=m$ (where $m,\, n\in\mathbb{N}$  and 
$m,\, n\geqslant 3$)  and $R_{1}R_{2}R_{3}$ is of finite order.
Then the candidates for $(n, m)$ will be $(3, 4)$, $(3, 5)$, $(4, 3)$, $(5, 4)$, $(8,6)$ 
and $(k, k)$ ($k\in\mathbb{N}$ and $k\geqslant3$).
\end{thm}

Note that the solutions correspond to the following parameter values, or their complex conjugates:
$$
\begin{array}{|l|lll|}
\hline
(n,m) & \rho & s=\rho-1 & \sigma=\tau  \\
\hline
(3,4) & (1+i\sqrt{7})/2 & e^{2\pi i/7}+e^{4\pi i/7}+e^{-6\pi i/7} 
& 1  \\
(3,5) & 2e^{2\pi i/5}\cos(\pi/5) & 
e^{2\pi i/5}+e^{7\pi i/15}+e^{-13\pi i/15} & 1  \\
(4,3) & 1 & 0 & \sqrt{2} \\
(5,4) & (1+i\sqrt{3})(\sqrt{5}-i\sqrt{3})/4 &
e^{-2\pi i/3}+e^{2\pi i/15}+e^{8\pi i/15} 
& (1+\sqrt{5})/2 \\
(8,6) & (1+i)(1-i/\sqrt{2}) &
e^{\pi i/2}+e^{\pi i/12}+e^{-7\pi i/12}
& \sqrt{2+\sqrt{2}} \\
(k,k) & 2e^{i\pi/k}\cos(\pi/k) & e^{2\pi i/k}  & 2\cos(\pi/k)  \\
\hline
\end{array}
$$

\section{The proof}
Firstly a direct computation will show that the symmetry $S$ conjugates $R_1R_2$ to $R_1R_3^{-1}R_2R_3$ and conjugates $R_2R_3$ to 
$R_1R_3$. It means that 
\[
br (R_1,\, R_2) = br (R_1,\, R_3^{-1}R_2R_3), \qquad br (R_2,\, R_3) = br (R_1,\, R_3)
\] 
by recalling Remark 2.2 and (\ref{equa:tr}). By the parameterization of the triangle groups in Section 2 and 
the assumption in Theorem \ref{thm:main}, one could get the matrix representation of $H,\,  R_1,\, R_2,\, R_3$ 
as (\ref{cal:H}), (\ref{mar:1}), (\ref{mar:2}), (\ref{mar:3}), where 
$$|\rho|=2\cos(\pi/m),\qquad \sigma=\tau=2\cos(\pi/n).$$
Throughout the proof we let $\zeta=\pi/m$ and $\eta=\pi/n$.

Because of $S^{2}=R_{1}R_{2}R_{3}$, we can restrict ourselves to $S$, which is elliptic of finite order. Equivalently, there exist $a$ and $b$ 
that are rational multiples of $\pi$ for which:
\begin{equation}\label{rel:1}
\tr(S)=-1+\rho=e^{i a}+e^{i b}+e^{-i (a+b)}.
\end{equation}
Observe that there is some ambiguity in the choice of $a$ and $b$. First, we can permute the
three terms in this expression, and so permute $\{a,\,b , \,-a-b\}$; secondly we can change the sign of
all three terms and, finally, since $\tr(S)$ is only defined up to multiplying by a cube root of unity,
we can add the same integer multiple of $2\pi/3$ to both $a$ and $b$. We will use these operations to
simplify things in our calculations below.

We denote $\tr(S)$ by $s$, then get that
\begin{eqnarray}
\label{rel:2}
|s|^2 = 1+|\rho|^2-2\Re(\rho) & = &|e^{i a}+e^{i b}+e^{-i (a+b)}|^2, \\
\label{rel:3}
\Re(s) = -1+\Re(\rho) & = & \cos(a)+\cos(b)+\cos(a+b),
\end{eqnarray}
Recall that 
\begin{eqnarray*}
|\rho|^2 & = & 4\cos^2\zeta=2\cos(2\zeta)+2, \\
\Re(\rho) & = & \frac{\sigma\tau}{2}=2\cos^2\eta=\cos(2\eta)+1. 
\end{eqnarray*}
The above two equations can be simplified to 
\begin{eqnarray}
\label{main}
1 & = & \cos(2\zeta)-\cos(2\eta)-\cos(a-b)-\cos(a+2b)-\cos(2a+b), \qquad \\
\label{minor}
0 & = & \cos(2\eta)-\cos(a)-\cos(b)-\cos(a+b).
\end{eqnarray}

In what follows we will repeatedly use the following result given by A. Monaghan, which generalises the result 
of Conway and Jones for vanishing sums of cosines of rational multiples of $\pi$.
 
\begin{prop}\text{\cite[Theorem 2.4.3.1]{Mon}}\label{posi:1} Suppose that we have at most five distinct 
rational  numbers of $\pi$, for which some rational linear combination of their cosines is rational but no proper 
subset has this property. If $\phi\in(0, \pi)$ and all other angles are normalised to lie in $(0, \frac{\pi}{2}),$ 
then the appropriate linear combination is proportional to one of the following:\\
\emph{(a)} $0=\sum_{j=0}^{2}\cos(\phi+\frac{2j\pi}{3})$,\\ 
\emph{(b)} $0=\cos(\phi)+\cos(\phi\pm\frac{2\pi}{5})-\cos(\phi\pm\frac{2\pi}{15})+\cos(\phi\pm\frac{7\pi}{15}),$\\
\emph{(c)} $0=\cos(\phi)-\cos(\phi\pm\frac{\pi}{5})+\cos(\phi\pm\frac{\pi}{15})-\cos(\phi\pm\frac{4\pi}{15}),$\\
\emph{(d)} $\frac{1}{2}=\cos(\frac{\pi}{3})$,\\
\emph{(e)} $\frac{1}{2}=\cos(\frac{\pi}{5})-\cos(\frac{2\pi}{5})$,\\
\emph{(f)} $\frac{1}{2}=\cos(\frac{\pi}{5})-\cos(\frac{\pi}{15})+\cos(\frac{4\pi}{15})$,\\
\emph{(g)} $\frac{1}{2}=-\cos(\frac{2\pi}{5})+\cos(\frac{2\pi}{15})-\cos(\frac{7\pi}{15})$,\\
\emph{(h)} $\frac{1}{2}=-\cos(\frac{\pi}{15})+\cos(\frac{2\pi}{15})+\cos(\frac{4\pi}{15})-\cos(\frac{7\pi}{15})$,\\
\emph{(i)}  $\frac{1}{2}=\cos(\frac{\pi}{7})-\cos(\frac{2\pi}{7})+\cos(\frac{3\pi}{7})$,\\
\emph{(j)}  $\frac{1}{2}=\cos(\frac{\pi}{7})-\cos(\frac{2\pi}{7})+\cos(\frac{2\pi}{21})-\cos(\frac{5\pi}{21})$,\\
\emph{(k)} $\frac{1}{2}=\cos(\frac{\pi}{7})+\cos(\frac{3\pi}{7})-\cos(\frac{\pi}{21})+\cos(\frac{8\pi}{21})$,\\
\emph{(l)} $\frac{1}{2}=-\cos(\frac{2\pi}{7})+\cos(\frac{3\pi}{7})+\cos(\frac{4\pi}{21})+\cos(\frac{10\pi}{21})$,\\
\emph{(m)} $\frac{1}{2}=\cos(\frac{\pi}{7})-\cos(\frac{\pi}{21})+\cos(\frac{2\pi}{21})-\cos(\frac{5\pi}{21})+\cos(\frac{8\pi}{21})$,\\
\emph{(n)} $\frac{1}{2}=-\cos(\frac{2\pi}{7})+\cos(\frac{2\pi}{21})+\cos(\frac{4\pi}{21})-\cos(\frac{5\pi}{21})+\cos(\frac{10\pi}{21})$,\\
\emph{(o)} $\frac{1}{2}=\cos(\frac{3\pi}{7})-\cos(\frac{\pi}{21})+\cos(\frac{4\pi}{21})+\cos(\frac{8\pi}{21})+\cos(\frac{10\pi}{21})$.\\
\end{prop}

 Since the right hand side of equation \eqref{main} is $1$ (rather than $0$ or $1/2$), Monaghan's 
 theorem implies that it must be a sum of (at least) two similar sums involving fewer cosines. We begin by
 showing that at least one of the cosines must itself be rational. 
 
 \begin{prop}\label{prop:one-cosine-rational}
 Suppose that $\zeta=\pi/m$, $\eta=\pi/n$ and $a$, $b$ are rational multiples of $\pi$ so that equations
 \eqref{main} and \eqref{minor} hold. Then one of the cosines in equation \eqref{main} must be rational.
 \end{prop}
 
 \begin{proof}
 Suppose that none of the cosines are rational. Then \eqref{main} splits into two rational sums, one of
 length two and the other of length three, neither of which has a rational subsum. 
 By inspection from prop~\ref{posi:1} we see that these two sums must have the value 
 $0$, $\pm 1/2$. Since they sum to $1$, they must both be $1/2$. Therefore, the sum of length 2 must be
 (e) and the sum of length 3 must be one of (f), (g) or (i).
  \begin{enumerate}
 \item[(1)] $1/2=\cos(2\zeta)-\cos(2\eta)=-\cos(a-b)-\cos(a+2b)-\cos(2a+b)$.
 Since $\zeta=\pi/m$ and $\eta=\pi/n$ the sum (e) implies that $2\zeta=\pi/5$ and $2\eta=2\pi/5$.
 For the second equation, there are certain symmetry operations on $a$ and $b$ 
 described in the paragraph after equation \eqref{rel:1} above.
Up to these operations, we now list the possible values of $a$ and $b$:
$$
\begin{array}{|c|c|c|c|c|}
\hline
a-b & a+2b & 2a+b & a & b \\
\hline
\pi/15 & 11\pi/15 & 4\pi/5 & 13\pi/45 & 2\pi/9 \\
2\pi/5 & 7\pi/15 & 13\pi/15 & 19\pi/45 & \pi/45 \\
2\pi/7 & 4\pi/7 & 6\pi/7 & -2\pi/7 & -4\pi/7 \\
\hline
\end{array}
$$
Using $2\eta=2\pi/5$, we see that none of the values in this table satisfy \eqref{minor}.
Therefore we get no solutions.
 \item[(2)] $1/2=\cos(2\zeta)-\cos(a-b)=-\cos(2\eta)-\cos(a+2b)-\cos(2a+b)$.
 The first equation gives $2\zeta=\pi/5$ as in case (1) and so $a-b=2\pi/5$.
 Since $a-b=(2a+b)-(a+2b)$ the difference of two of the angles in the second equation
 must be $2\pi/5$. By inspection, we see the only solution is 
 $a+2b=7\pi/15$ and $2a+b=13\pi/15$. This means $2\eta=2\pi/5$ and we are
 back in case (1). 
 \item[(3)] $1/2=-\cos(2\eta)-\cos(a-b)=\cos(2\zeta)-\cos(a+2b)-\cos(2a+b)$.
 The first equation gives $2\eta=2\pi/5$ as in (1) and so $a-b=4\pi/5$. Substituting in
 the second equation, we see $a+2b=-\pi/15$ and $2a+b=11\pi/15$. Thus 
 $2\zeta=\pi/5$ and we are back in case (1) again.
 \item[(4)] $1/2=-\cos(a-b)-\cos(a+2b)=\cos(2\zeta)-\cos(2\eta)-\cos(2a+b)$.
 Up to symmetries of $a$, $b$ and $-a-b$, the first sum implies that 
 $a-b=2\pi/5$ and $a+2b=4\pi/5$. Hence $2a+b=6\pi/5$ and so the second sum must
 be (f). Thus $\cos(2\zeta)=\cos(2\pi/m)=\cos(4\pi/15)$ or 
 $\cos(2\eta)=\cos(2\pi/n)=-\cos(4\pi/15)$ so either $m$ or $n$ is not an integer.
 Therefore there are no solutions.
\end{enumerate}
 \end{proof}
 
 \medskip
 
 As a consequence of this result, we can consider separate cases where each of the cosines
 in \eqref{main} is rational. If either $\cos(2\zeta)$ or $\cos(2\eta)$ is rational it must be $0$ or 
 $\pm 1/2$ since $\zeta=\pi/m$ and $\eta=\pi/n$ where $m$ and $n$ are at least 3. If one of 
 the other three cosines is rational we can use the allowable symmetries of $a$ and $b$, we 
 to assume that $\cos(a-b)$ is rational. We treat each of these cases separately below. 
 First we eliminate a simple situation which gives us many solutions and will recur in the different
 cases.
 
 \begin{lem}\label{lem-n=m}
Suppose that  $\cos(2\zeta)=\cos(2\eta)$, or equivalently $m=n$, then putting $s=e^{\pm 2\pi i/m}$ 
gives a solution to equations \eqref{main} and \eqref{minor} for all $m\geqslant 3$. 
\end{lem}

\begin{proof}
Substituting $\cos(2\zeta)=\cos(2\eta)$ into \eqref{main} gives:
Now we obtain that
\begin{eqnarray*}
0 & = & 1+\cos(a-b)+\cos(a+2b)+\cos(2a+b) \\
& = & 2\cos^2\bigl((a-b)/2\bigr)+2\cos\bigl((a-b)/2\bigr)\cos\bigl(3(a+b)/2\bigr) \\
& = & 4\cos\bigl((a-b)/2\bigr)\cos\bigl((a+2b)/2\bigr)\cos\bigl((2a+b)/2\bigr).
\end{eqnarray*}
Therefore one of $(a-b)$, $(a+2b)$ or $(2a+b)$ is an odd multiple of $\pi$.
Without loss of generality, we suppose that $a+2b=(2k+1)\pi$.
Then we get  $-a-b=b-(2k+1)\pi$ which yields  $s=e^{ia},$ where $a$ 
is a rational multiple of $\pi$. Because $\Re(s)=-1+\Re(\rho)=-1+\frac{|\sigma|^2}{2}=\cos(2\eta)$, 
we see that $\cos(a)=\cos(2\pi/m)=\cos(2\pi/n)$.

Now we consider the signature of the Hermitian form 
\[
\Det(H)=ie^{-\frac{4\theta+3\phi}{2}i}(-1+e^{(2\theta+\phi)i})(e^{i\theta}+e^{i\phi})^2.
\]
\begin{center}
\begin{table}[h]
\caption{Signature of Hermitian form}\label{tab:3}
\begin{tabular}{|c|c|c|c|}
\hline
$s$ & (2, 1) & degenerate & (3, 0) \\
\hline
$e^{\frac{2 \pi i}{3}}$ ($m=n=3$) &$p\geqslant4$& $p=3$ & $p=2$\\\hline
$e^{\frac{2 \pi i}{4}}$ ($m=n=4$) &$p\geqslant3$& $p=2$ & none \\\hline
$e^{\frac{2 \pi i}{k}}$ ($m=n=k\geqslant5$)&$p\geqslant2$& none & none\\\hline
\end{tabular}
\end{table}
\end{center}
In this case, we get the solution $n=m$.\\
\end{proof}

\medskip

We now consider the cases where $\cos(2\zeta)$, $\cos(2\eta)$ or $\cos(a-b)$ are rational.
We will use the following result proved by Parker when he was analysing the triangle groups
with 3-fold symmetry \cite{Par3}. In \cite{Par3} the last two cases were missed out, but this was corrected
in \cite{DPP2}.

\begin{prop}\text{\cite[prop 3.2]{Par3}}\label{posi:2}
Let $\theta$, $a$ and $b$ be rational multiples of $\pi$. Write  $s=e^{ia}+e^{ib}+e^{-i(a+b)}$. 
Then the only possible solutions to the equation
$$
\cos(2\theta)-\cos(a-b)-\cos(a+2b)-\cos(2a+b)=\frac{1}{2}
$$
give rise to the following values of $\theta$ and $s$, up to changing the sign of $\theta$
and up to conjugating $s$ and multiplying it by a power of $\omega=e^{2\pi i/3}$:
\begin{enumerate}
\item[(\romannumeral1)] $2\theta=2\pi/3$ and  $s=-e^{-i \psi/3}$ for some angle $\psi$ that is a 
rational multiple of $\pi$;
\item[(\romannumeral2)] $2\theta=\psi$ and  $s=e^{2i\psi/3}+e^{-i\psi/3}=e^{i\psi/6}2\cos\frac{\psi}{2}$ 
for some angle $\psi$ that is a rational multiple of $\pi$;
\item[(\romannumeral3)] $2\theta=\pi/3$ and $s=e^{i\pi/3}+e^{-i\pi/6}2\cos\frac{\pi}{4}$;
\item[(\romannumeral4)] $2\theta=\pi/5$ and $s=e^{i\pi/3}+e^{-i\pi/6}2\cos\frac{\pi}{5}$;
\item[(\romannumeral5)] $2\theta=3\pi/5$ and $s=e^{i\pi/3}+e^{-i\pi/6}2\cos\frac{2\pi}{5}$;
\item[(\romannumeral6)] $2\theta=\pi/2$ and $s=e^{2\pi i/7}+e^{4\pi i/7}+e^{-6 \pi i/7}$;
\item[(\romannumeral7)] $2\theta=\pi/2$ and $s=e^{2\pi i/9}+e^{-i \pi/9}2\cos\frac{2\pi}{5}$;
\item[(\romannumeral8)] $2\theta=\pi/2$ and $s=e^{2\pi i/9}+e^{-i \pi/9}2\cos\frac{4\pi}{5}$;
\item[(\romannumeral9)] $2\theta=\pi/7$ and $s=e^{2\pi i/9}+e^{-i \pi/9}2\cos\frac{2\pi}{7}$;
\item[(\romannumeral10)] $2\theta=5\pi/7$ and $s=e^{2\pi i/9}+e^{-i \pi/9}2\cos\frac{4\pi}{7}$;
\item[(\romannumeral11)] $2\theta=3\pi/7$ and $s=e^{2\pi i/9}+e^{-i \pi/9}2\cos\frac{6\pi}{7}$;
\item[(\romannumeral12)] $2\theta=2\pi/5$ and $s=1+2\cos\frac{2\pi}{5}$;
\item[(\romannumeral13)] $2\theta=4\pi/5$ and $s=1+2\cos\frac{4\pi}{5}$.
\end{enumerate}
\end{prop}

Note that for the groups Parker was considering $s=e^{ia}+e^{ib}+e^{-ia-ib}$ was the trace of $R_1J$,
whereas in our case it is the trace of $S$. In the cases where $\cos(2\zeta)=1/2$ or 
$\cos(2\eta)=-1/2$ then equation \eqref{main} reduces to the equation from prop~\ref{posi:2},
and we can use that result to find solutions.

\begin{lem}\label{lem-cos3eta-rational}
Suppose that $\cos(2\eta)$ is rational. Then the only solutions to \eqref{main} and \eqref{minor}
are $\cos(2\zeta)=\cos(2\pi/m)$ and $\cos(2\eta)=\cos(2\pi/n)$ where
$(n,m)$ is one of $(3,3)$, $(3,4)$, $(3,5)$, $(4,3)$, $(4,4)$ or $(6,6)$.
\end{lem}

\begin{proof}
Since $\cos(2\eta)$ is rational and not equal to $\pm1$ it can only be $0$ or $\pm1/2$. We treat each case
separately. 
\begin{enumerate}
\item[(1)]
$\cos(2\eta)=-\frac{1}{2}$, which gives $n=3.$ Note that 
$$
\frac{1}{2}=\cos(2\eta)+1=\Re(\rho)=\Re(s)+1
$$
and so $\Re(s)=-1/2$.
We rewrite (\ref{main}) to give the equation from prop~\ref{posi:2} with $\theta=\zeta$.
%
%
%
%
By direct calculation, we just need to consider cases (\romannumeral1), (\romannumeral2) and 
(\romannumeral6) because of $\Re(s)=-1/2$.
\begin{enumerate}
\item[(\romannumeral1)] $s=-e^{-i\psi/3}$ and so $|s|=1$. This yields that $|\rho|=2\cos(\pi/m)=1$, 
and so $m=3$. By considering 
$\Re(s)=-\cos(\theta/3)$, we know that $\theta=\pm\pi+6k\pi$ ($k\in\mathbb{Z}$) which means that 
$s=-e^{\mp i \pi/3}$. From (\ref{H:verify}), we get that
\[
\Det(H)=\mp \sqrt{3}\cos(\phi/2)+\sin(\phi/2)-2\sin(3\phi/2).
\]
We list the corresponding signature of Hermitian form for different $s$ in Table 3.2. 
\begin{center}
\begin{table}[h]\label{tab:(1)}
\caption{Signature of Hermitian form}
\begin{tabular}{|c|c|c|c|}
\hline
$s$ & (2, 1) & degenerate & (3, 0) \\
\hline
$-e^{-i\pi/3}$&$p\geqslant4$& $p=3$& $p=2$\\\hline
$-e^{i\pi/3}$&none& $p=6$& $p\neq6$\\\hline
\end{tabular}
\end{table}
\end{center}
In this case, we get that $n=m=3$.
\item[\noindent(\romannumeral2) ]
$s=e^{2i\psi/3}+e^{-i\psi/3}=e^{i\psi/6}2\cos(\psi/2)$ where $\psi=2\theta$.
By solving
\begin{eqnarray*}
-1/2 & = & \Re(s)=\cos(4\theta/3)+\cos(2\theta/3) \\
& = & 2\cos^2(2\theta/3)+\cos(2\theta/3)-1
\end{eqnarray*}
we obtain $\cos(2\theta/3)=(-1\pm\sqrt{5})/4$. That is, $2\theta/3=\pm 2\pi/5+2k\pi$ or 
$\pm4\pi/5+2k\pi$. Hence $\theta=\pm3\pi/5+3k\pi$ or $\pm6\pi/5+3k\pi$.
The only solution to $|s|=2|\cos(\theta)|=2\cos(\pi/m)$ is $m=5$ 
(coming from $2\theta/3=4\pi/5-2\pi$).
Therefore $s=e^{2\pi i/5}+e^{4\pi i/5}$ or $s=e^{-2\pi i/5}+e^{-4\pi i/5}$.
%
In these cases we find, respectively, that:
\begin{align*}
\Det(H)&=-\sqrt{5+2\sqrt{5}}\cos{\frac{\phi}{2}}-(2+\sqrt{5}+4\cos{\phi})\sin{\frac{\phi}{2}},
\\
\Det(H)&=\sqrt{5+2\sqrt{5}}\cos{\frac{\phi}{2}}-(2+\sqrt{5}+4\cos{\phi})\sin{\frac{\phi}{2}}.
\end{align*}

\begin{center}
\begin{table}[h]
\caption{Signature of Hermitian form}
\begin{tabular}{|c|c|c|c|}
\hline
$s$ & (2, 1) & degenerate & (3, 0) \\
\hline
$e^{-\frac{2\pi i}{5}}+e^{-\frac{4\pi i}{5}}$&$p\leqslant7$& none & $p\geqslant8$\\\hline
$e^{\frac{2\pi i}{5}}+e^{\frac{4\pi i}{5}}$&$p\geqslant2$& none& none\\\hline
\end{tabular}
\end{table}
\end{center}
In this case, we get that $n=3,$ $m=5$.

\item[\noindent(\romannumeral6)] 
$s=e^{2\pi i/7}+e^{4\pi i/7}+e^{-6\pi i/7}=(-1+\sqrt{7}i)/2$. It follows that $\Re(s)=-1/2$ and 
$|s|=\sqrt{2}$ which indicates that $m=4.$ A simple calculation yields that
\[
\Det(H)=\frac{1}{2}(1-8\cos\phi)\sin\frac{\phi}{2}
\]
from which it follows that the signature of the Hermitian form will be of $(2, 1)$ for $p\geqslant5$, otherwise 
it will be positive. In this case, we get that $n=3$, $m=4$.
\end{enumerate}
Therefore we obtain the solutions $(n,m)=(3,3)$, $(3,4)$ and $(3,5)$.
\item[(2)] $\cos(2\eta)=0.$ Now we have $|\sigma|^2=2$ which yields $\Re(s)=0$. 
Therefore one can get the following two equations
\begin{equation}
\left\{
\begin{aligned}
&\cos(2\zeta)-\cos(a-b)-\cos(a+2b)-\cos(2a+b)=1,\\
&\cos a+\cos b+\cos(a+b)=0.
\end{aligned}
\right.
\end{equation}
Since the first of these has 1 on the right hand side, it must split as the sum of (at least) two minimal
subsums. Treating these case by case we see that the only possibilities are 
$\cos(2\zeta)=0$, which yields $m=n=4$ and $\cos(2\zeta)=-1/2$, which gives $n=4$ and $m=3$.
The former case is a particular instance of lem~\ref{lem-n=m}.
In the latter case
we rewrite \eqref{rel:2} as
$$
|s|^2=1+|\rho|^2-2\Re(\rho)=1+2\cos(2\zeta)-2\cos(2\eta)=0.
$$
Therefore, the only solution is $s=0$, or equivalently $\rho=1$.
This implies that
\[
\Det(H)=-2\sin\frac{3\phi}{2}=-2\sin\frac{3\pi}{p},
\]
and the signature of the Hermitian form will be positive if $p=2$, degenerate if 
$p=3$, negative (of signature (2,1)) if $p\geqslant4.$ Therefore in this case we get that $n=4$, $m=3$.
Hence the only solutions we get in this case are $(n,m)=(4,3)$ and $(4,4)$.
\item[(3)] $\cos(2\eta)=1/2.$ Now we have 
$|\sigma|^2=3$ from which it follows that $\Re(s)=1/2$. We rewrite the two equations
\begin{equation}
\left\{
\begin{aligned}
&\cos(2\zeta)-\cos(a-b)-\cos(a+2b)-\cos(2a+b)=\frac{3}{2},\\
&\cos a+\cos b+\cos(a+b)=\frac{1}{2}.
\end{aligned}
\right.
\end{equation}
If the second equation is irreducible, then it must be one of prop~\ref{posi:1} parts (f), (g) or (i).
We see in each case that the angles involved do not sum to $0$ (making each cosine positive,
the sum is $\pi$ times the ratio of two odd integers for each choice of sign). 
If the second equation splits as the sum of two rational subsums then, without loss of generality,
$\cos(a)$ is rational. Hence it is in the set $\{0,\pm1/2.\pm1\}$. Simple trigonometry shows that
\begin{eqnarray*}
2\cos(a/2)\cos(a/2+b)& = & \cos(b)+\cos(a+b)\\
& = & 1/2-\cos(a),\\
\cos(a+2b)+1 & = & 2\cos^2(a/2+b) \\
& = & \bigl(1/2-\cos(a)\bigr)^2/\bigl(1+\cos(a)\bigr), \\
\cos(a-b)+\cos(2a+b) & = & 2\cos(3a/2)\cos(a/2+b)\\ 
& = & -2\bigl(1/2-\cos(a)\bigr)^2.
\end{eqnarray*}
Substituting these identities in the first equation, we see that 
$\cos(2\zeta)$ is a rational function of $\cos(a)$, and so is rational. Substituting the different values
of $\cos(a)$ gives a solution with $\zeta=\pi/m$ only when $\cos(a)=\pm1/2$. In both cases,
$\cos(2\zeta)=1/2$ and so $m=6$.
Thus we obtain the solution $(n,m)=(6,6)$.
\end{enumerate}
\end{proof}

\medskip

\begin{lem}\label{lem-cos2zeta-rational}
Suppose that $\cos(2\zeta)$ is rational. Then the only solutions to \eqref{main} and \eqref{minor}
are $\cos(2\zeta)=\cos(2\pi/m)$ and $\cos(2\eta)=\cos(2\pi/n)$ where
$(n,m)$ is one of $(3,3)$, $(4,3)$, $(4,4)$, $(5,4)$, $(6,6)$ or $(8,6)$.
\end{lem}

\begin{proof}
Since $\cos(2\zeta)$ is rational and not equal to $\pm1$ it can only be $0$ or $\pm1/2$. We treat each case
separately. 
\begin{enumerate}
\item[(1)]
$\cos(2\zeta)=1/2$, which gives $m=6$. In this case, we know
$$
|s|^2=1+2\cos(2\zeta)-2\cos(2\eta)=2-2\cos(2\eta).
$$
In this case, we rewrite equation \eqref{main} to give the equation from prop~\ref{posi:2}
with $2\theta=\pi-2\eta$.
Checking one by one, we will find that there is no value of $s$ in prop \ref{posi:2} satisfying 
(\ref{minor}) except the cases (\romannumeral1) and (\romannumeral2). For
(\romannumeral1) we have $2\eta=\pi-2\theta=\pi/3$ and so $n=6$ (we have analysed this case previously). 
For (\romannumeral2), we have $\psi=\pi-2\eta$ and 
$s
=e^{2i\pi/3-4i\eta/3}+e^{-i\pi/3+2i\eta/3}$.
Substituting in equation \eqref{minor}
gives
\begin{eqnarray*}
0 & = & \cos(2\eta)-{\rm Re}(s) \\
& = & -\cos(\pi-2\eta)-\cos(2\pi/3-4\eta/3)-\cos(\pi/3-2\eta/3) \\
& = & -\cos(2\pi/3-4\eta/3)\bigl(1+2\cos(\pi/3-2\eta/3)\bigr).
\end{eqnarray*}
The only solution with $\eta=\pi/n$ is when $2\pi/3-4\eta/3=\pi/2$. That is, $n=8$. By calculating $\Det(H)=-2\cos(\phi)(1+2\sin\phi)$, we see that $H$ is of signature $(3, 0)$ for $p=2$ and is of signature $(2, 1)$ for any $p\geqslant3$.
In this case, we get $(n,m)=(6,6)$ or $(8,6)$.\\
\item[(2)] $\cos(2\zeta)=0$, which gives $m=4$.
Then we get that $|\rho|^2=2$ and $|s|^2=3-|\sigma|^2$. Also 
(\ref{main}) can be replaced by
\[
-\cos(2\eta)-\cos(a-b)-\cos(a+2b)-\cos(2a+b)=1.
\]
We have already analysed the case where $\cos(2\eta)=0$ or $-1/2$, which lead to the solution
$(n,m)=(3,4)$ or $(4,4)$. If $\cos(2\eta)=1/2$, 
then $|\sigma|^2=3$ induces $s=0$, which contradicts $\Re(s)=-1+|\sigma|^2/2=1/2$. 
Then it suffices for us to 
consider the following possible values due to $\eta=\pi/n$,
\begin{center}
\begin{tabular}{|c|c|c|c|c|c|}
\hline
$2\eta$ & $a-b$ & $a+2b$ & $2a+b$ & $a$ & $b$ \\
\hline
$2\pi/5$ &$2\pi/3$& $7\pi/15$ & $17\pi/15$ & $3\pi/5$& $-\pi/15$\\\hline
$2\pi/5$ &$2\pi/5$& $4\pi/5$ & $6\pi/5$ & $8\pi/5$& $2\pi/15$\\\hline
\end{tabular}
\end{center}
From this table, we know that $n=5$ and the pair values $a=3\pi/5$ and $b=-\pi/15$ do not satisfy the equation (\ref{minor}). However  the second line $a=8\pi/5$ and $b=2\pi/15$ satisfy the equation (\ref{minor}) by applying the equation 
(g) in prop \ref{posi:1}. Then we calculate the signature of the Hermitian form $H$ using (\ref{H:verify}). We see that $H$ is of signature $(3, 0)$ for $p=2$ and is of signature $(2, 1)$ for any $p\geqslant 3$.
In this case, we get that $m=4$ and $n=5$.\\

\item[(3)] $\cos(2\zeta)=-\frac{1}{2}.$ It follows that $m=3$ and
\[
-\cos(2\eta)-\cos(a-b)-\cos(a+2b)-\cos(2a+b)=\frac{3}{2}.
\]
Also, $\cos(2\zeta)=-1/2$ implies $|\rho|=1$ and so $\cos(2\eta)+1=\Re(\rho)\leqslant 1$.
This means that $\cos(2\eta)\le 0$ and so either $\cos(2\eta)=\cos(2\pi/n)=-1/2$ or $0$.
We have analysed both of these cases already. These give solutions $(n,m)=(3,3)$ or $(4,3)$.
\end{enumerate}
\end{proof}

\medskip
Now we begin to consider the remaining case in which $\cos(a-b)$ is rational. 
\begin{lem}\label{lem-1}
Suppose that $\cos(a-b)=-1$, then $\cos(2\zeta)-\cos(2\eta)=0$, and the possible solutions are given in lem~\ref{lem-n=m} in which $n=m$.
\end{lem}
\begin{proof}
$\cos(a-b)=-1$ which gives 
$b=a+(2k+1)\pi$.  Hence we have $\cos(a+2b)=\cos(3a)$ and $\cos(2a+b)=-\cos(3a)$.
Therefore, equation (\ref{main}) reduces to 
$\cos(2\zeta)-\cos(2\eta)=0$, which we have already treated in 
lem~\ref{lem-n=m}.
\end{proof}
\begin{lem}\label{lem-2}
Suppose that $\cos(a-b)=-1/2$, $\cos(2\zeta)$ and $\cos(2\eta)$ are not rational, 
$\cos(2\zeta)-\cos(2\eta)\neq0$. 
Then we get no solutions for $n$, $m$ such that (\ref{main}) and (\ref{minor}) hold.
\end{lem}
\begin{proof}
$\cos(a-b)=-1/2$, which gives
$b=a\pm 2\pi/3+2k\pi$. Hence we have $\cos(a+2b)=\cos(3a\mp2\pi/3)$ and $\cos(2a+b)=\cos(3a\pm2\pi/3)$.
Therefore equation (\ref{main}) becomes
\begin{eqnarray*}
1/2 & = & 1+\cos(a-b) \\
& = & \cos(2\zeta)-\cos(2\eta)-\cos(a+2b)-\cos(2a+b) \\
& = & \cos(2\zeta)-\cos(2\eta)+\cos(3a).
\end{eqnarray*}
Since we have supposed that $\cos(2\zeta)$ and $\cos(2\eta)$ are not rational, the only way
this equation can split into to rational subsums is for $\cos(3a)$ to be rational. Investigating the
different possibilities, we see that \eqref{minor} then implies $\cos(2\eta)$ is rational.

Now suppose the equation does not split into two rational sums of cosines.
We list the possible values of $2\zeta,\ 2\eta,\ a,\, b$, 
up to the allowable symmetries of $a$ and $b$. 
\begin{center}
\begin{table}[h]
\begin{tabular}{|c|c|c|c|c|c|}
\hline
$2\zeta$ & $2\eta$ & $3a$ & $a-b$ & $a$ & $b$ \\
\hline
$\pi/5$ &$2\pi/5$& $\pi/2$& $2\pi/3$& $\pi/6$ & $-\pi/2$ \\\hline
$2\pi/5$ &$\pi/5$& $0$& $2\pi/3$& $0$ & $-2\pi/3$ \\\hline
$\pi/5$ &$\pi/15$& $4\pi/15$& $2\pi/3$& $4\pi/45$ & $-26\pi/45$ \\\hline
$2\pi/15$ &$2\pi/5$& $8\pi/15$ & $2\pi/3$& $8\pi/45$ & $-22\pi/45$ \\\hline
$\pi/7$ &$2\pi/7$& $3\pi/7$ & $2\pi/3$& $\pi/7$ & $-11\pi/21$ \\\hline
\end{tabular}
\end{table}
\end{center}
However, we note that there are no values of $2\eta,\ a,\, b$ in this list satisfying (\ref{minor}). Therefore there are no solutions for $n, m$.
\end{proof}

\begin{lem}\label{lem-3}
Suppose that $\cos(a-b)=0$,  $1/2$ or $1$, $\cos(2\zeta)$ and $\cos(2\eta)$ are not rational, $\cos(2\zeta)-\cos(2\eta)\neq0$. Then there are no solutions for $n$, $m$ satisfying both (\ref{main}) and (\ref{minor}).
\end{lem}
\begin{proof}
We immediately get that
\begin{equation}\label{eq:cos(a-b)}
\cos(2\zeta)-\cos(2\eta)-\cos(a+2b)-\cos(2a+b)=1\hbox{ or }\frac{3}{2}\hbox{ or }2.
\end{equation}
Since the right hand side is not $0$, $\pm 1/2$, we see that this sum must split into shorter rational
sums of cosines.We break down into the following three cases.
\begin{enumerate}
\item[(1)] $\cos(2\zeta)-\cos(2\eta)=\pm1/2$.
 \begin{enumerate}
 \item[(\romannumeral1)] $\cos(2\zeta)-\cos(2\eta)=1/2$. Note that $\zeta=\pi/m$ and $\eta=\pi/n$, where $m, n\in\mathbb{N}$. Therefore we know that $(n,m)$ is  $(5,10)$ and
\[
\cos(a-b)+\cos(a+2b)+\cos(2a+b)=-\frac{1}{2}.
\]
We have supposed that $\cos(a-b)$ is rational, then using
elementary trigonometry arguments, we see that 
\[
2\cos\bigl((a-b)/2\bigr)\cos\bigl(3(a+b)/2\bigr)=-\frac{1}{2}-\cos(a-b).
\]
Squaring both sides and rearranging gives
$$
\cos(3a+3b)=\frac{\cos^2(a-b)-3/4}{\cos(a-b)+1}.
$$
We have assumed that $\cos(a-b)=0$ or  $\cos(a-b)=1/2$ or $\cos(a-b)=1$, which means that $\cos(3a+3b)=-3/4$ or $-1/3$ or $-1/8$. It gives a contradiction here.\\
 \item[(\romannumeral2)] $\cos(2\zeta)-\cos(2\eta)=-1/2$, i.e. 
 $$
 \cos(a+2b)-\cos(2a+b)=-\frac{3}{2}\hbox{ or }-2\hbox{ or }-\frac{5}{2}.
 $$
 This sum must again split and so both cosines are rational. 
Therefore the possible values for $\cos(a+2b)$ are just $-1$ or $-1/2$ which are equivalent to the case where $\cos(a-b)$ is this value, see lem~\ref{lem-1} and lem \ref{lem-2}. However we have assumed that $\cos(2\zeta)-\cos(2\eta)\neq0$, therefore there are no solutions for $n$, $m$ satisfying both (\ref{main}) and (\ref{minor}).
\end{enumerate}
\item[(2)] $\cos(2\zeta)-\cos(x)$ (or $\cos(2\eta)+\cos(y)$) 
is $1/2$ or $-1/2$, where $x,\,y\in\{a+2b,\,2a+b\}$. 

Recalling the equation (\ref{eq:cos(a-b)}), $\cos(2\zeta)-\cos(x)=\pm1/2$ means that $\cos(2\eta)+\cos(y)$ is one of the values $\{-5/2,~-2,~-3/2,~-1,~-1/2\}$. We just need to consider the case $\cos(2\eta)+\cos(y)=-1/2$, because other values of $\cos(2\eta)+\cos(y)$ mean that $\cos(2\eta)$ will be rational. Without loss of generality, we suppose that $x=a+2b$, $y=2a+b$ and list the values of $2\zeta$, $a+2b$, $2\eta$, $2a+b$ and corresponding $a-b$ .
\begin{center}
\begin{table}[h]
\begin{tabular}{|c|c|c|c|c|}
\hline
$2\zeta$ & $a+2b$ & $2\eta$ & $2a+b$ & $a-b$ \\
\hline
%
$\pi/5$ &$2\pi/5$& $2\pi/5$ & $4\pi/5$& $2\pi/5$ \\\hline
$\pi/5$ &$-2\pi/5$& $2\pi/5$ & $4\pi/5$& $6\pi/5$ \\\hline
\end{tabular}
\end{table}
\end{center}
There are no values of $a-b$ such that $\cos(a-b)=0$ or  $\cos(a-b)=1/2$ or $\cos(a-b)=1$.\\
Therefore there are no solutions for $n$, $m$ satisfying both (\ref{main}) and (\ref{minor}) in this case.

\item[(3)] $\cos(x)$ is rational, where $x\in\{a+2b,\, 2a+b\}$. By suitable changes of $a$ and $b$, the cases $\cos(a+2b)$ or $\cos(2a+b)$ is $-1/2$ or $-1$
are equivalent to the cases in lem \ref{lem-1} and lem \ref{lem-2}. Therefore there are no solutions for $n$, $m$ because we supposed that $\cos(2\zeta)-\cos(2\eta)\neq0$. 

Then we consider the condition for $\cos(x)$ to be $0$, $1/2$ or $1$ and suppose that $x=a+2b.$ We get that 
$$\cos(2\zeta)-\cos(2\eta)-\cos(2a+b)\in\left\{1,~\frac{3}{2},~2,~\frac{5}{2},~3\right\},$$
which can be reduced to $\cos(2\zeta)-\cos(2a+b)$ or $\cos(2\eta)+\cos(2a+b)$ is rational which has been considered above.
\end{enumerate}
Now we can get that there are no solutions for $n$, $m$ satisfying both (\ref{main}) and (\ref{minor}) under the conditions in lem \ref{lem-3}.
\end{proof}

\noindent We sum up all the candidates for $n$, $m$ from above process,

lem 3.3 \quad $n=m\geqslant3$;

lem 3.5 \quad $(n,m)\in (3,3)$, $(3,4)$, $(3,5)$, $(4,3)$, $(4,4)$ or $(6,6)$;

lem 3.6 \quad $(n,m)\in (3,3)$,  $(4,3)$, $(4,4)$, $(5,4)$, $(6,6)$ or $(8,6)$;\\
which we desired.

\begin{rem}
Note that the new candidates for $(n, m)$ to be $(5, 4)$, $(4, 3)$ and $(8,6)$ do not appear on
Thompson's list in \cite{Thom}. However referring to \cite{DPP2}, in what follows we will see that 
the triangle groups for $(n, m)$ to be $(5, 4)$ corresponds to Thompson groups ${\bf S}_2$ and 
the triangle groups for $(m, n)$ to be $(4, 3)$ is of actually Mostow groups with braiding $(2, 3, 4; 4)$.
The pair $(n,m)=(8,6)$ was also found by Deraux when he was making a similar computer search to
Thompson (private communication). 
\vskip 0.5\baselineskip
\noindent1. $(n, m)=(5, 4)$.\\
Suppose that $M_1, \, M_2,\, M_3$ are three complex reflections of order $p$, which satisfy 
\[
br( M_1, M_2)=4,\quad br(M_1, M_3)=br(M_2, M_3)=3,\quad br(M_1, M_3^{-1}M_2M_3)=5.
\] 
Actually, $M_1, \, M_2,\, M_3$ will be Thompson group ${\bf S}_2$. Write 
$R_1=M_2^{-1}M_1M_2$, $R_2=M_1M_2M_1^{-1}$, $R_3=M_3$. We claim that
\[
br(R_1, R_2)=br(R_1, R_3^{-1}R_2R_3)=4, \quad br(R_1, R_3)=br(R_2, R_3)=5.
\]
First, observe, we also have $br(M_2^{-1}M_1M_2, M_3)=br(M_1^{-1}M_2M_1, M_3)=5$. Thus
\[
br(R_1, R_3)=br(M_2^{-1}M_1M_2, M_3)=5,\quad br(R_2, R_3)=br(M_1M_2M_1^{-1}, M_3)=5.
\]
Using $br(M_1, M_2)=4$, we have
\[
R_1R_2=(M_2^{-1}M_1M_2)(M_1M_2M_1^{-1})=M_2^{-1}(M_2M_1M_2M_1)M_1^{-1}=M_1M_2.
\]
Hence $br(R_1, R_2)=br(M_1, M_2)=4.$ We denote $M_1$, $M_1^{-1}$ by $1$, $\bar{1}$ simply 
and so on. Now we consider 
\begin{align*}
R_1R_3^{-1}R_2R_3
&= M_2^{-1}M_1M_2M_3^{-1}M_1M_2M_1^{-1}M_3\\
&= \bar{2}12\bar{3}12\bar{1}3\\
&= (123123)(\bar{3}\bar{2}\bar{1}\bar{3}\bar{2}\bar{1} \cdot 
\bar{2}12\bar{3}12\bar{1}3 \cdot 123123)\bar{3}\bar{2}\bar{1}\bar{3}\bar{2}\bar{1}\\
&= (123123)(\bar{3}\bar{2}\bar{1}\bar{3} \cdot 12\bar{1} \cdot 
\bar{3}12\bar{1}3 \cdot 123123)\bar{3}\bar{2}\bar{1}\bar{3}\bar{2}\bar{1}\\
&= (123123)(\bar{3}\bar{2}(\bar{1}\bar{3}1)2(\bar{1}\bar{3}1)2
(\bar{1}31)23123)\bar{3}\bar{2}\bar{1}\bar{3}\bar{2}\bar{1}\\
&= (123123)(\bar{3}\bar{1}31\bar{2}\bar{1}\bar{3}13123)\bar{3}\bar{2}\bar{1}\bar{3}\bar{2}\bar{1}\\
&= (123123)(1\bar{3}\bar{2}323)\bar{3}\bar{2}\bar{1}\bar{3}\bar{2}\bar{1}\\
&= (123123)(12)\bar{3}\bar{2}\bar{1}\bar{3}\bar{2}\bar{1}\\
\end{align*}
Since $R_1R_3^{-1}R_2R_3$ is conjugate to $M_1M_2$ we see that 
$br(R_1, R_3^{-1}R_2R_3)=br(M_1, M_2)=4$ as claimed.
In particular this shows that  this case is equivalent to Thompson groups ${\bf S}_2$.
\vskip 0.5\baselineskip
\noindent2. $(n, m)=(4, 3)$\\
In this case, it is easy to check that $R_{2}$, $R_{3}$ (also $R_{3}$, $R_{1}$) braid with length 
4, $R_1$, $R_2$ (also $R_1$, $R_3^{-1} R_2 R_3$) braid with length 3, 
$R_1$, $R_2 R_3 R_2^{-1}$ (also $R_3$, $R_1 R_2 R_1^{-1}$) braid with length 2 
(i.e. they commute)
and $R_1 R_2 R_3$ is regular elliptic of order 3. Note that $\Det(H)<0$ when $p\geqslant4$.

As the same fashion in \cite{KPT1}, we define $\iota$ by the reflection of group that acts on the 
generating set ($R_{1},\, R_{2},\, R_{3}$) as follows,
$$\iota(R_{1})=R_{1},\quad \iota(R_{2})=R_{1}R_{2}R_{1}^{-1}, \quad \iota(R_{3})=R_{3}.$$
Under the action of $\iota$, the (4, 4, 3; 3)-triangle groups will be sent to the triangle groups with 
braiding (2, 3, 4; 4) 
\begin{align*}
\langle\iota(R_{1}),\, \iota(R_{2})\, \iota(R_{3}): &\iota(R_2 R_3)=\iota(R_3 R_2),
(\iota(R_1 R_2))^{\frac{3}{2}}=(\iota(R_2 R_1))^{\frac{3}{2}},\\
&(\iota(R_1 R_3))^2=(\iota(R_3 R_1))^2,\\
&(\iota(R_{1}R_{2}R_{3}R_{2}^{-1}))^2= (\iota(R_{2}R_{3}R_{2}^{-1}R_{1}))^2\rangle.
\end{align*}

Recall the Mostow groups $\Gamma(p, t)$ mentioned in \cite{Mos, Par3}. For Mostow groups, 
there exists a complex hyperbolic isometry $J$ of order 3 so that $R_{j+1}=JR_{j}J^{-1}$ and 
$R_{i}R_{i+1}R_{i}=R_{i+1}R_{i}R_{i+1}$. We could rewrite them as triangle groups with 
braiding (2, 3, 4; 4) as follows
\begin{align*}
\langle R_{1}, R_{2}, J(R_{1}R_{2})^{-1}:&
R_2 J(R_1R_{2})^{-1}=J(R_1R_{2})^{-1} R_2,
(R_1 R_2)^{\frac{3}{2}}=(R_2 R_1)^{\frac{3}{2}},\\
&(R_1 J(R_1R_{2})^{-1})^2=(J(R_1R_{2})^{-1} R_1)^2, \\
&(R_1 R_2 J(R_1 R_2)^{-1} R_2^{-1})^2=(R_2 J(R_1 R_2)^{-1} R_2^{-1} R_1)^2
\rangle.
\end{align*}
\end{rem}

\noindent\textbf{Acknowledgement} This research is carried out while John Parker was visiting Department of Mathematical and Computing Sciences 
at Tokyo Institute of Technology. The first author is supported by a JSPS Invitation Fellowship. The second author is 
supported by Grant-in-Aid Scientific Research (B) No. 15H03619 (Sadayoshi Kojima, Tokyo Institute of Technology) 
and Grant-in-Aid Scientific Research (S) No. 24224002 (Takashi Tsuboi, University of Tokyo). We also thank Martin Deraux
for several useful conversations.

\bibliographystyle{siam}
\bibliography{Sun-bib.bib}

\end{document}